\documentclass[11pt]{amsart}

\usepackage[square,compress,comma, numbers,sort]{natbib}
\usepackage[colorlinks=true, citecolor=red, linkcolor=blue]{hyperref}
\usepackage{amsfonts,mathtools}

\allowdisplaybreaks[4]

\usepackage{amsmath}
\usepackage{bbm}
\usepackage{amssymb}
\usepackage{mathtools}

\usepackage{color}
\definecolor{c20}{rgb}{0.,0.7,0.}
\definecolor{c30}{rgb}{0.,0.,1.}
\definecolor{c40}{rgb}{1,0.1,0.7}
\definecolor{c50}{rgb}{1,0,0}
\definecolor{c60}{rgb}{1,0.9,0.1}

\def\Var{\text{Var}}

\definecolor{c20}{rgb}{0.,0.7,0.}
\definecolor{c30}{rgb}{0.,0.,1.}
\definecolor{c40}{rgb}{1,0.1,0.7}
\definecolor{c50}{rgb}{1,0,0}
\definecolor{c60}{rgb}{1,0.9,0.1}

\newcommand{\ve}{\varepsilon}

\newcommand{\E}[1]{\mathbb{E}\left\{ #1\right\}}

\newcommand{\pk}[1]{\mathbb{P} \left\{ #1 \right \} }
\newcommand{\PP}{\mathbb{P}}

\newcommand{\R}{\mathbb{R}}

\newcommand{\inr}{\in \R}

\newcommand{\mD}{\mathbb{D}}

\newcommand{\limit}[1]{\lim_{#1 \to   \infty}}

\newcommand{\BQN}{\begin{eqnarray}}
\newcommand{\EQN}{\end{eqnarray}}
\newcommand{\BQNY}{\begin{eqnarray*}}
\newcommand{\EQNY}{\end{eqnarray*}}

\newcommand{\BS}{\begin{sat}}
\newcommand{\ES}{\end{sat}}
\newcommand{\BT}{\begin{theo}}
\newcommand{\ET}{\end{theo}}
\newcommand{\BK}{\begin{korr}}
\newcommand{\EK}{\end{korr}}

\DeclareMathOperator{\cov}{cov}

\newcommand{\BD}{\begin{de}}
\newcommand{\ED}{\end{de}}

\newcommand{\BIT}{\begin{itemize}}
\newcommand{\EIT}{\end{itemize}}
\newcommand{\BDI}{\begin{description}}
\newcommand{\EDI}{\end{description}}

\newcommand{\BRM}{\begin{remarks}}
\newcommand{\ERM}{\end{remarks}}

\newcommand{\BEL}{\begin{lem}}
\newcommand{\EEL}{\end{lem}}

\def\bqny#1{{\begin{eqnarray*} #1 \end{eqnarray*}}}
\def\bqn#1{{\begin{eqnarray} #1 \end{eqnarray}}}

\newtheorem{theo}{Theorem}[section]
\newtheorem{sat}[theo]{Proposition}
\newtheorem{de}[theo]{Definition}
\newtheorem{lem}[theo]{Lemma}

\newtheorem{korr}[theo]{Corollary}

\newtheorem{remarks}[theo]{Remarks}
\newtheorem{prop}[theo]{Proposition}

\newcommand{\COM}[1]{}

\def\IF{\infty}

\newcommand{\QED}{\hfill $\Box$}

\def\ve{\varepsilon}

\def\IF{\infty}

\newcommand{\pr}{\mathcal{P}_{T_u}}
\newcommand{\Cr}{\mathbb{C}_{T_u}}
\newcommand{\ind}{\mathbb{I}}

\begin{document}

\COM{
We would like to point out seminal manuscript \cite{DI2005}
establishing asymptotics of $\lambda(u)$
under week assumptions on variance and covariance of $X$.
For the
discrete-time investigations (i.e., when $t$ in model \eqref{*Risk_Model}
belongs to a discrete grid $\{0,\delta,2\delta,...\}$ for some $\delta>0$),
we refer to \cite{HashorvaGrigori,Piterbargdiscrete,
secondproj}.
We would like to suggest a reader contributions
\cite{KrzysPeng2015,DebickiScaledBM,KEP20171,
ParisianFiniteHoryzon,SumFBMDebicki,
Rolski17,Bai_2017,DHJ15}
for the related generalizations of the classical ruin problem
and monographs \cite{20lectures,Pit96,MR1379077} for
a survey of the known results.
}

\title{Sojourn Ruin of a Two-Dimensional
Fractional Bronwian Motion Risk Process}

\author{{Grigori Jasnovidov}}
\address{Grigori Jasnovidov, Department of Actuarial Science,
University of Lausanne,\\
UNIL-Dorigny, 1015 Lausanne, Switzerland
}
\email{griga1995@yandex.ru}

 \maketitle

{\bf Abstract:} This paper derives the asymptotic behavior of
$$\pk{\int\limits_0^\IF \ind\Big(B_H(s)-c_1s>q_1u,
B_H(s)-c_2s>q_2u\Big)ds>T_u},\quad u \to \IF,$$
where $B_H$ is a fractional Brownian motion,
$c_1,c_2,q_1,q_2>0,\  H \in (0,1), \ T_u \ge 0$ is a measurable function and
$\ind(\cdot)$ is the indicator function.

\\


\section{Introduction \& Preliminaries}

Consider the risk model defined by
\bqn{\label{*Risk_Model}
R(t) = u+\rho t-X(t), \ \ \ \ t \ge 0,
}
where $X(t)$ is a centered Gaussian risk process
with a.s. continuous sample paths,
$\rho>0$ is the net profit rate and $u>0$ is the initial capital.
This model is relevant to insurance and financial applications, see,
e.g., \cite{MR1458613}.
A question of numerous investigations
(see \cite{DI2005,HashorvaGrigori,Piterbargdiscrete,
secondproj,KrzysPeng2015,DebickiScaledBM,KEP20171,
ParisianFiniteHoryzon,SumFBMDebicki,
Rolski17,Bai_2017,DHJ15,20lectures,Pit96,MR1379077,PiterbargFilomat})
 is the study of the asymptotics of the classical ruin
probability
\bqn{\label{*classical_ruin_probability}
\lambda(u) := \pk{\exists t\ge 0: R(t)<0}
}
as $u \to \IF$ under different levels of generality.
It turns out, that only
for $X$ being a Brownian motion (later on BM) $\lambda(u)$ can be calculated explicitly:
if $X$ is a standard BM, then $\lambda(u) = e ^{-2\rho u}, \ u,\rho>0,$
see \cite{DeM15}. Since it seems impossible to find the exact value
of $\lambda(u)$ in other cases,
the approximations of $\lambda(u)$ as $u \to \IF$
is dealt with.
Some contributions (see, e.g., \cite{SojournHashorvaCorrelatedBM,
SojournInfty}), extend the
classical ruin problem to the so-called sojourn ruin problem, i.e.,
approximation of the sojourn ruin probability defined by
\bqn{\label{*paris_ruin-def}
\pk{\int\limits_0^\IF
\ind (R(s)<0)ds>T_u},
}
where $T_u\ge 0$ is a measurable function of $u$.
As in the classical case, only for $X$ being a BM
 the probability above can
be calculated explicitly, see \cite{SojournInfty}:
$$\pk{\int\limits_0^\IF
\ind (B(s)-cs>u)ds>T}=
 \Big(2(1+c^2T)\Psi(c\sqrt T)-\frac{c\sqrt{2T}}{\sqrt \pi}
 e^{\frac{-c^2T}{2}}\Big)
e^{-2cu}
, \ \ c>0, \ T,u\ge 0,$$
where $\Psi$ is the survival function of a standard
Gaussian random variable, $B$ is a standard BM and $\ind(\cdot)$ is
the indicator function. Motivated by
\cite{Lanpeng2BM} (see also \cite{secondproj,thirdprojectParisian}),
we study a generalization of the main problem in \cite{Lanpeng2BM}
for the sojourn ruin, i.e., we shall study the asymptotics of

$$\Cr(u) :=
\pk{\int\limits_0^\IF  \ind\Big(
B_H(s)-c_1s>q_1u,B_H(s)-c_2s>q_2u\Big)ds>T_u},
 \ \ \ \  u \to \IF,$$
where $B_H$ is a standard fractional Brownian motion (later on fBM), i.e.,
a Gaussian process with zero expectation and covariance defined by
\bqny{
\cov(B_H(s),B_H(t)) = \frac{|t|^{2H}+|s|^{2H}-|t-s|^{2H}}{2},
 \ \ \ \ \ \ t,s \inr.
}
The ruin probability above is of interest for
reinsurance models, see \cite{Lanpeng2BM} and references therein.
By the self-similarity of fBM we have
\bqny{
\Cr(u) &=&
\pk{\int\limits_0^\IF\ind\Big(
B_H(su)>c_1su+q_1u,B_H(su)>c_2su+q_2u\Big)d(su)>T_u}
\\&=&
\pk{\int\limits_0^\IF\ind\Big( u^{H}B_H(s)
>(c_1s+q_1)u,u^HB_H(s)>(c_2s+q_2)u\Big)ds>T_u/u}
\\&=&
\pk{\int\limits_0^\IF \ind\left(\frac{B_H(s)}
{\max(c_1s+q_1,c_2s+q_2)}> u^{1-H}\right)ds>T_u/u}
.}

 In order to prevent the problem of degenerating
to the one-dimensional sojourn problem
discussed in \cite{SojournHashorvaCorrelatedBM,SojournInfty}
(i.e., to impose the denominator in the line above be nonlinear function)
we assume that
\bqn{\label{*cq}
c_1>c_2, \ \ \ q_2>q_1.
}
The variance of
the process two lines above can achieve its unique maxima only
at one of the following points:
\bqn{\label{*ttdefinition}
t_* = \frac{q_2-q_1}{c_1-c_2}, \ \ \ t_1 = \frac{q_1H}{(1-H)c_1}, \ \ \
t_2 = \frac{q_2H}{(1-H)c_2}.
}

From \eqref{*cq} it follows that $t_1<t_2$;
as we shall see later, the order between $ t_1,t_2$ and $t_*$
determines the asymptotics of $\Cr(u)$ as $u \to \IF$.
As mentioned in \cite{ParisianFiniteHoryzon},
for the approximation of the one-dimensional
Parisian ruin probability we need to control
the growth of $T_u$ as $u \to \IF$. As in \cite{ParisianFiniteHoryzon},
we impose the following condition:
\bqn{\label{*assumption_T_u}
 \limit{u}T_u u^{1/H-2} =T \in [0,\IF), \ H\in (0,1)
.}
Note that $T_u$ satisfying \eqref{*assumption_T_u}
may go to $\IF$ for $H>1/2$, converges to non-negative limit for $H=1/2$
and approaches $0$ for $H<1/2$ as $u\to\IF$.
We see later on in Proposition \ref{*remark_H<1/2} that the condition
above is necessary and it seems very difficult to derive
the exact asymptotics of $\Cr(u)$ as $u\to \IF$ without it.
\newline

The rest of the paper is organized in the following way. In the next
section we present the main results of the paper, in Section
\ref{section_proofs} we give all proofs, while technical
calculations are deferred to the Appendix.

\section{Main Results}

Define for some function $h$ and $K\ge 0$ the sojourn Piterbarg constant by
\bqny{
\mathcal{B}^h_K
=
\int\limits_\R
\pk{\int\limits_{-\IF}^{\IF}\ind\Big(\sqrt 2 B(s)-|s|+h(s)>x\Big)ds
 >K}e^{x}dx
}
when the integral above is finite and Berman's constant by
$$
\mathcal B_{2H}(x) =
\limit{S}\frac{1}{S}\int\limits_{\R} \pk{
\int\limits_0^S\ind(\sqrt 2 B_H(t)-t^{2H}+z>0)dt>x
}e^{-z}dz, \quad x\ge 0.$$
It is known (see, e.g., \cite{SojournInfty}) that
$\mathcal B_{2H}(x) \in (0,\IF)$ for all $x\ge 0$;
we refer to \cite{SojournInfty}
and references therein for the properties of relevant Berman's
constants.
Let for $i=1,2$
\bqn{\label{*mathbb_D}
 \ \ \mD_H = \frac{c_1t_*+q_1}{t_*^H}, \
K_H = \frac{2^{\frac{1}{2}-\frac{1}{2H}}\sqrt \pi}{\sqrt{H(1-H)}},
\
\mathbb{C}_{H}^{(i)} = \frac{c_i^Hq_i^{1-H}}{H^H(1-H)^{1-H}},
\ D_i =
\frac{c_i^2(1-H)^{2-\frac{1}{H}}}{2^{\frac{1}{2H}}H^{2}}.
}

Now we are ready to give the asymptotics of $\Cr(u)$ as $u\to \IF$.
\begin{theo} \label{*theoparis}
 Assume that \eqref{*cq} holds and $T_u$ satisfies \eqref{*assumption_T_u}.\\
1)If $t_*\notin (t_1,t_2)$, then as $u \to \IF$
\bqn{\label{*parisclaim11}
\mathbb{C}_{T_u}(u) \sim (\frac{1}{2})^{\mathbb{I}(t_*=t_i)}\times
\begin{cases}
\Big(2(1+c_i^2T)\Psi(c_i\sqrt T)-\frac{c_i\sqrt{2T}}{\sqrt \pi}
e^{\frac{-c_i^2T}{2}}\Big)
e^{-2c_iq_iu},& H=1/2\\
K_H\mathcal{B}_{2H}(TD_i)(\mathbb{C}_H^{(i)}u^{1-H})^{\frac{1}{H}-1}
\Psi(\mathbb{C}_H^{(i)}u^{1-H}),
 & H\neq 1/2,
\end{cases}}
where $i=1$ if $t_*\le t_1$ and $i=2$ if $t_*\ge t_2$.

2) If $t_*\in (t_1,t_2)$ and $\lim\limits_{u\to \IF}T_uu^{2-1/H} = 0$ for
$H>1/2$, then as $u\to\IF$
\bqn{\label{*parisclaim2}
\Cr(u) \sim \Psi(\mD_H u^{1-H})\times
\begin{cases}
1,&H>1/2\\
\mathcal{B}_{T'}^d,& H=1/2\\
\mathcal{B}_{2H}(\overline D T)
 Au^{(1-H)(1/H-2)},& H<1/2,
\end{cases}
}
where $\mathcal{B}_{T'}^d \in (0,\IF)$,
\bqn{\label{*defT'd}
T' = T\frac{(c_1q_2-q_1c_2)^2}{2(c_1-c_2)^2}, \
\ d(s) =
s\frac{c_1q_2+c_2q_1-2c_2q_2}{c_1q_2-q_1c_2}\mathbb{I}(s<0)
+s\frac{2c_1q_1-c_1q_2-q_1c_2}{c_1q_2-q_1c_2}\mathbb{I}(s\ge0)
}
 and
\bqn{\label{*A_overline_D}\!\!
A = \Big(|H(c_1t_*+q_1)-c_1t_*|^{-1} +
|H(c_2t_*+q_2)-c_2t_*|^{-1}\Big) \frac{t_*^H\mD_H^{\frac{1}{H}-1}}
{2^{\frac{1}{2H}}}, \ \ \
\overline D=\frac{(c_1t_*+q_1)^{\frac{1}{H}}}{2^{
\frac{1}{2H}
}t_*^2}. }
\end{theo}
Note that if $T=0$, then
the result above reduces to Theorem 3.1 in \cite{Lanpeng2BM}.
As already mentioned in the introduction \eqref{*assumption_T_u}
is a necessary condition for the theorem above.
To illustrate situation when it is not satisfied
we consider a "simple" scenario with $T_u$ being a positive
constant.
\begin{prop} \label{*remark_H<1/2}
If $H<1/2, \ T_u=T>0$ and $t_* \in (t_1,t_2)$, then
\bqny{
\bar{C} \Psi(\mD_H u^{1-H})e^{-C_{1,\alpha} u^{2-4H}-
 C_{2,\alpha}u^{2(1-3H)}}
&\le& \Cr(u)\notag
\\&\le&
(2+o(1))\Psi(\mD_Hu^{1-H})\Psi\left(u^{1-2H}\frac{T^H\mD_H}{2t_*^H}\right)
,\quad u \to \IF,}
where $\bar C\in (0,1)$ is a fixed constant that does not depend on $u$ and
\bqny{
\alpha = \frac{T^{2H}}{2t_*^{2H}}, \ \ \ \ \ \
C_{i,\alpha} = \frac{\alpha^{i}}{i} \mD_H^2, \ \ \
i=1,2
.}
\end{prop}
Note that the lower bound in the proposition above
decays to zero exponentially faster than the upper bound as $u\to \IF$.

\section{Proofs}\label{section_proofs}
First we give the following auxiliary results.
As shown, e.g., in Lemma 2.1 in \cite{PickandsA}
\bqn{\label{ratio_uniform}
(1-\frac{1}{u^2})
\frac{1}{\sqrt{2\pi} u}e^{-u^2/2}
\le \Psi(u) \le \frac{1}{\sqrt{2\pi} u}e^{-u^2/2}, \ \ \ u >0.
}

Recall that $K_H,D_1$ and $\mathbb{C}_H^{(1)}$ are defined in
\eqref{*mathbb_D}.
A proof of the proposition below is given in the Appendix.
\begin{prop}\label{*prop_sojourn_app}
Assume that $T_u$ satisfies \eqref{*assumption_T_u}.
Then as $u\to \IF$
\bqny{\label{*theo_app}
 \pk{\int\limits_ 0^\IF \ind(B_H(t)-c_1t>q_1u)dt>T_u}
\sim
\begin{cases}
\Big(2(1+c_1^2T)\Psi(c_1\sqrt T)-\frac{c_1\sqrt{2T}}{\sqrt \pi}
 e^{\frac{-c_1^2T}{2}}\Big)
e^{-2c_1q_1u},
&H=1/2,\\
K_H\mathcal{B}_{2H}(TD_1)(\mathbb{C}_H^{(1)}u^{1-H})^{\frac{1}{H}-1}
\Psi(\mathbb{C}_H^{(1)}u^{1-H}),
& H\neq 1/2.
\end{cases}}
\end{prop}

Now we are ready to perform our proofs. \\

\textbf{Proof of Theorem \ref{*theoparis}.}
\textbf{Case (1).} Assume that $t_*<t_1$.
Let
\bqny{
V_i(t) = \frac{B_H(t)}{c_it+q_i}
\ \ \ \text{ and } \ \ \
\psi_i(T_u,u) = \pk{\int\limits_0^\IF
\ind(B_H(t)-c_it>q_iu)ds>T_u } , \ \ \ i =1,2.
}
For $0<\ve<t_1-t_*$ by the self-similarity of fBM
we have
\bqny{
\psi_1(T_u,u) \ge \Cr(u)
\ge
\pk{\int\limits_{t_1-\ve}^{t_1+\ve} \ind( V_1(t)> u^{1-H},
V_2(t)> u^{1-H})dt>\frac{T_u}{u}}
=
\pk{\int\limits_{t_1-\ve}^{t_1+\ve} \ind( V_1(t)> u^{1-H})dt>\frac{T_u}{u}}.
}
We have by Borel-TIS inequality, see \cite{20lectures}
(details are in the Appendix)
\bqn{\label{*mainintpariscase1}
\psi_1(T_u,u) \sim
\pk{\int\limits_{t_1-\ve}^{t_1+\ve} \ind(
V_1(t)> u^{1-H})ds>T_u/u }, \ \ \ u \to \IF
}
implying $\Cr(u)\sim \psi_1(T_u,u)$ as $u \to \IF$.
The asymptotics of $\psi_1(T_u,u)$ is given
in Proposition \ref{*prop_sojourn_app}, thus the claim follows.
\\
\\
Assume that $t_*=t_1$.
We have
\bqny{
\pk{\int\limits_{t_1}^\IF
\ind(V_1(s)> u^{1-H})ds>T_u}
&\le& \Cr(u)
\\&\le&
\pk{\int\limits_{t_1}^\IF
\ind(V_1(s)> u^{1-H})ds>T_u}+
\pk{\!\exists t\!\in \![0,t_1]:
V_2(t)\!>\! u^{1-H}}
.}
From the proof of Theorem 3.1,
case (4) in \cite{Lanpeng2BM} it follows that
the second term in the last line above
is negligible comparing with the final asymptotics of
$\Cr(u)$ given in \eqref{*parisclaim11}, hence
$$\Cr(u) \sim \pk{\int\limits_{t_1}^\IF
\ind(V_1(s)> u^{1-H})ds>T_u}, \quad u \to \IF.$$
Since $t_1$ is the unique maxima of $\Var\{V_1(t)\}$
from the proof of Theorem 2.1, case i) in \cite{SojournInfty} we have
\bqny{
\pk{\int\limits_{t_1}^\IF
\ind(V_1(t)>u^{1-H})dt>T_u/u}
&\sim& \frac{1}{2}\pk{\int\limits_{0}^\IF
\ind(V_1(t)>u^{1-H})dt>T_u/u}
\\&=&
\frac{1}{2}\pk{\int\limits_0^\IF
\ind(B_H(t)-c_1t>q_1u)dt>T_u}
, \ \ u \to \IF.}
The asymptotics of the last probability above is given in Proposition
\ref{*prop_sojourn_app} establishing the claim.
Case $t_*\ge t_2$ follows by the same arguments.
\\

\textbf{Case (2).}
Assume that $H>1/2$.
We have by Theorem 2.1 in \cite{thirdprojectParisian} and Theorem 3.1 in \cite{Lanpeng2BM}
with \bqny{
\mathcal{R}_{T_u}(u) &=& \pk{\exists t\ge 0: B_H(t)-c_1t>q_1u,
B_H(t)-c_2t>q_2u},\\
\pr(u) &=& \pk{\exists t\ge 0: \inf\limits_{s\in [t,t+T_u]}
(B_H(s)-c_1s)>q_1u, \inf\limits_{s\in [t,t+T_u]}(B_H(s)-c_2 s)>q_2 u} }
that
$$\Psi(\mD_H u^{1-H}) \sim
\pr(u)\le \Cr(u)\le \mathcal{R}_{T_u}(u)\sim
\Psi(\mD_H u^{1-H}), \quad u \to \IF,$$
and the claim follows.
\\
\\
Assume that $H=1/2$. First let \eqref{*assumption_T_u} holds with
$T_u=T>0$. We have as $u \to \IF$ and then $S \to \IF$ (proof is in the
Appendix)
\bqn{\label{*infintcum}
\Cr(u) \sim \pk{\int\limits_{ut_*-S}^{ut_*+S}
\ind\Big(B(s)-c_1s>q_1u,B(s)-c_2s>q_2u\Big)ds>T}
=: \kappa_S(u)
.}
Next with $\phi_u$ the density of $B(ut_*)$, $
\eta = c_1t_*+q_1=c_2t_*+q_2$ and $\eta_*=\eta/t_*-c_2=q_2/t_*$ we have
\bqny{
\kappa_S(u)
&=&
\int\limits_\R
\pk{
\int\limits_{ut_*-S}^{ut_*}\ind(B(s)-c_2s>q_2u)ds+\!\!\!\!
\int\limits_{ut_*}^{ut_*+S}\!\!\! \ind(B(s)-c_1s>q_1u)ds>T
\big|B(ut_*)=\eta u-x} \phi_u(\eta u -x)dx
\\&=&
\int\limits_\R \PP\Big\{
\int\limits_{ut_*-S}^{ut_*} \ind (B(s)-c_2s>q_2u)ds
\\ & \ & +
\int\limits_{ut_*}^{ut_*+S}\!\!\! \ind(B(s)-B(ut_*)-c_1(s-ut_*)
-c_1ut_*
>q_1u-\eta u+x)ds>T\big
|B(ut_*)=\eta u-x\Big\}
\phi_u(\eta u -x)dx
\\ &=&
\int\limits_\R
\pk{
\int\limits_{ut_*-S}^{ut_*}\ind(B(s)-c_2s>q_2u)ds+
\int\limits_0^S \ind(B_*(s)-c_1s>x)ds>T\big
|B(ut_*)=\eta u-x}
\phi_u(\eta u -x)dx
\\ &=&
\frac{e^{-\frac{\eta^2 u}{2t_*}}}{\sqrt{2\pi u t_*}}
\int\limits_\R
\pk{\int\limits_{-S}^{0}\ind\Big(Z_u(s)+\eta_*s>x\Big)ds+
\int\limits_0^S \ind(B_*(s)-c_1s>x)ds >T}
e^{\frac{\eta x}{t_*}-\frac{x^2}{2ut_*}}dx
,}
where $Z_u(t)$ is a Gaussian process with expectation and covariance
defined below:
\bqn{\label{*Z_u}
\E{Z_u(t)} = \frac{-x}{ut_*}t, \ \ \ \ \ \
\cov (Z_u(s),Z_u(t)) = \frac{-st}{ut_*}-t, \ \ \ \ \  s \leq t \le 0
.}
Since $Z_u(t)$ converges to BM in the sense of convergence
finite-dimensional
distributions for any fixed $x \inr$ as $u \to \IF$ we have (details are in the Appendix)
\bqn{\label{*cumasym_elimin_u}
& \ &
\int\limits_\R
\pk{\int\limits_{-S}^{0}\ind\Big(Z_u(s)+\eta_*s>x\Big)ds+
\int\limits_0^S \ind(B_*(s)-c_1s>x)ds >T}e^{\frac{\eta x}{t_*}-
\frac{x^2}{2ut_*}}dx
\notag\\ &\sim&
\int\limits_\R
\pk{\int\limits_{-S}^{0}\ind\Big(B(s)+\eta_*s>x\Big)ds+
\int\limits_0^S \ind(B_*(s)-c_1s>x)ds >T}e^{\frac{\eta x}{t_*}}dx
\\&=:&\notag K(S).
}
Since $\pk{\exists t\ge 0: B(t)-ct>x} = e^{-2cx}, \ c,x>0$ (see, e.g.,
 \cite{DeM15}) we have
\bqny{
K(S) &\le&
\int\limits_0^\IF \Big(\pk{\exists s<0: B(s)+\eta_*s>x}+
\pk{\exists s\ge 0 : B_*(s)-c_1s>x}\Big) e^{\frac{\eta x}{t_*}}dx
+ \int\limits_{-\IF}^0 e^{\frac{\eta x}{t_*}}dx
\\&=& \int\limits_0^\IF \Big(
e^{(-2\eta_*+\eta/t_*)x}+e^{(-2c_1+\eta/t_*)x}\Big)dx
+t_*/\eta < \IF
}
provided by $t_*\in (t_1,t_2)$. Since $K(S)$ is an increasing function and
$\lim\limits_{S \to \IF}K(S) < \IF$ we have as $S \to \IF$
\bqny{
K(S) & \to &
\int\limits_\R
\pk{\int\limits_{0}^{\IF}\ind\Big(B(s)-\eta_*s>x\Big)ds+
\int\limits_0^\IF \ind(B_*(s)-c_1s>x)ds >T}e^{\frac{\eta x}{t_*}}dx
\\&=&
\frac{t_*}{\eta}
\int\limits_\R
\pk{\int\limits_{0}^{\IF}\ind\Big(B(s)-\frac{\eta_*t_*}{\eta}s>x\Big)ds+
\int\limits_0^\IF \ind(B_*(s)-\frac{c_1t_*}{\eta}s>x)ds
>\frac{\eta^2T}{t_*^2}}e^xdx
\notag\\&=&
\frac{t_*}{\eta}\int\limits_\R
\pk{\int\limits_{-\IF}^{\IF}\ind\Big(\sqrt 2B(s)-|s|+d(s)>x\Big)ds
>\frac{\eta^2T}{2t_*^2}}e^xdx
\notag\\&=&
\frac{t_*}{\eta}\mathcal{B}_{T'}^d \in (0,\IF)
\notag,}
where $T'$ and $d(s)$ are defined in \eqref{*defT'd}.
Finally, combining \eqref{*cumasym_elimin_u} with
the line above we have as $u \to \IF$ and then $S \to \IF$
$$\kappa_S(u) \sim
\mathcal{B}_{T'}^d
\Psi(\mathbb{D}_{1/2}\sqrt u)$$
and by \eqref{*infintcum} the claim follows.
If \eqref{*assumption_T_u} holds with $T_u = 0$, then we obtain the claim
immediately by Theorem 3.1 in \cite{Lanpeng2BM} and observation that
$\mathcal{B}_{0}^d$ coincides with the corresponding Piterbarg constant
introduced in \cite{Lanpeng2BM}.
\\

Now assume that \eqref{*assumption_T_u} holds with
any possible $T_u$. If \eqref{*assumption_T_u} holds with $T>0$, then for
large $u$ and any $\ve>0$ it holds, that
$\mathbb{C}_{(1+\ve)T}(u)\le \Cr(u) \le
\mathbb{C}_{(1-\ve)T}(u)$ and hence
$$(1+o(1))\mathcal{B}_{ T'(1+\ve)}^d\Psi(\mathbb{D}_{1/2}\sqrt u)\le
\Cr(u)\le \mathcal{B}_{ T'(1-\ve)}^d
\Psi(\mathbb{D}_{1/2}\sqrt u)(1+o(1)),\qquad u \to \IF.$$
By Lemma 4.1 in \cite{SojournInfty}
 $\mathcal{B}_{x}^d$ is a continuous function with respect
to $x$ and thus letting $\ve \to 0$ we obtain the claim.
If \eqref{*assumption_T_u} holds with $T=0$, then for
large $u$ and any $\ve>0$ we have
$$\mathcal{B}_{\ve}^d\Psi(\mathbb{D}_{1/2}\sqrt u)\le
\Cr(u)\le \mathcal{B}_{0}^d\Psi(\mathbb{D}_{1/2}\sqrt u)$$
and again letting $\ve \to 0$ we obtain the claim by continuity of
$\mathcal{B}_{(\cdot)}^d$.
\newline

Assume that $H<1/2$.
First we have with $\delta_u = u^{2H-2}\ln^2 u$ as
$u \to \IF$ (proof is in Appendix)
\bqn{\label{*cum_proof_main_int_H<1/2}
\notag\Cr(u) &\sim&
\pk{\int\limits_{ut_*-u\delta_u}^{ut_*}
\ind(B_H(t)-c_2t>q_2u)dt>T_u }
+ \pk{\int\limits_{ut_*}^{ut_*+u\delta_u}
\ind(B_H(t)-c_1t>q_1u)dt>T_u }
\\&=:& g_1(u)+g_2(u)
.}

Assume that \eqref{*assumption_T_u} holds with $T>0$.
Using the approach from \cite{SojournInfty}
we have with $\ind_a(b) = \ind(b>a), \ a,b\inr$
\bqny{
g_2(u)&=& \pk{\int\limits_0^{\delta_u T_u^{-1}u}\ind_{M(u)}
\Big(\frac{B_H(ut_*+tT_u)}{u(q_1+c_1t_*)+c_1tT_u}M(u) \Big)dt>1}
\\&=:& \pk{ \int\limits_0^{\delta_u T_u^{-1}u}\ind_{M(u)}(Z_u^{(1)}(t))dt>1}
\\&=&
\pk{ \int\limits_0^{\delta_u T_u^{-1}u
K_1}
\ind_{M(u)}(Z_u^{(1)}(tK_1^{-1}))dt>K_1}
\\&=:&
\pk{ \int\limits_0^{\delta_u T_u^{-1}u
K_1}
\ind_{M(u)}(Z_u^{(2)}(t))dt>K_1},
}
where
$$
K_1=\frac{T\mD_H^{1/H}}{2^{\frac{1}{2H}}t_*}, \quad
M(u) = \inf_{t\in [t_*,\IF)}\frac{u(c_1t+q_1)}{\Var\{B_H(ut)\}}
 = \mD_Hu^{1-H}.
$$
For variance $\sigma_{Z_u^{(2)}}^2(t)$ and correlation $r_{Z_u^{(2)}}(s,t)$
 of $Z_u^{(2)}$
 for $t,s \in [0,\delta_u T_u^{-1}uK_1]$ it holds, that as $u\to\IF$
\bqny{
1-\sigma_{Z_u^{(2)}}(t) &=& \frac{2^{\frac{1}{2H}}t_*^H\mD_H^{1-1/H}
|q_1H-(1-H)c_1t_*|}
{(q_1+c_1t_*)^2}tu^{1-1/H}+O(t^2u^{2(1-1/H)}), \\
 1-r_{Z_u^{(2)}}(s,t) &=& \mD_H^{-2}u^{2H-2}|t-s|^{2H}+O(u^{2H-2}|t-s|^{2H}
 \delta_u).
}
Now we apply Theorem 2.1 in \cite{SojournInfty}. All conditions of
the theorem are fulfilled with parameters
\bqny{ & \ &
\omega(x) = x, \ \overleftarrow{\omega}(x) = x,
 \ \beta =1, \ g(u) =
\frac{2^{\frac{1}{2H}}t_*^H\mD_H^{1-1/H}
|q_1H-(1-H)c_1t_*|}
{(q_1+c_1t_*)^2}u^{1-1/H},
\\& \ &
\eta_\varphi(t) = B_H(t),\ \sigma_\eta^2(t) = t^{2H}, \
\Delta(u)=1, \ \varphi = 1,
\\& \ &
 n(u) = \mD_H u^{1-H}, \
 a_1(u) = 0, \ a_2(u) = \delta_u T_u^{-1}uK_1, \ \gamma = 0, \
x_1 = 0, \ x_2 = \infty, \ y_1=0, \ y_2 = \IF, \ x = K_1,
\\& \ &
\theta(u) = u^{(1/H-2)(1-H)} \mD_H^{-1+1/H}
|q_1H-(1-H)c_1t_*|^{-1}t_*^{H}2^{-\frac{1}{2H}},
}
and thus as $u \to \IF$
$$
g_2(u) = \pk{ \!\!\int\limits_0^{\delta_u T_u^{-1}u
K_1}\!\!\!\!\!\!\!\!
\ind_{M(u)}(Z_u^{(2)}(t))dt>K_1} \sim
\mathcal{B}_{2H}(\frac{T\mD_H^{\frac{1}{H}}}
{2^{\frac{1}{2H}}t_*})u^{(\frac{1}{H}-2)(1-H)}
\frac{t_*^{H}\mD_H^{-1+1/H}}{2^{\frac{1}{2H}}|q_1H-(1-H)c_1t_*|}
\Psi(\mD_Hu^{1-H}).
$$
Similarly we obtain
$$
g_1(u) \sim
\mathcal{B}_{2H}(\frac{T\mD_H^{1/H}}
{2^{\frac{1}{2H}}t_*})u^{(1/H-2)(1-H)}
\frac{t_*^{H}\mD_H^{-1+1/H}}{2^{\frac{1}{2H}}|q_2H-(1-H)c_2t_*|}
\Psi(\mD_Hu^{1-H}), \ \ \ u \to \IF
$$
and the claim follows if in \eqref{*assumption_T_u} $T>0$.
 Now let
$\eqref{*assumption_T_u}$ holds with $T=0$.
Since $\pr(u)\le\Cr(u)\le \mathcal{R}_{T_u}(u)$
we obtain the claim by Theorem 2.1 in \cite{thirdprojectParisian}
 and Theorem 3.1 in \cite{Lanpeng2BM}. \QED
\\

\textbf{Proof of Proposition \ref{*remark_H<1/2}.}
The proof of this proposition is the same as the proof of Proposition
2.2 in \cite{thirdprojectParisian}, thus we refer to
\cite{thirdprojectParisian} for the proof. \QED

\section{Appendix}

\textbf{Proof of \eqref{*mainintpariscase1}.}
To establish the claim we need to show that
\bqny{
\pk{\int\limits_{[0,\IF)\backslash [t_1-\ve,t_1+\ve]}
\ind\left(V_1(s)> u^{1-H}\right)ds>T_u/u} = o(\psi_1(T_u,u)), \ \
u \to \IF.
}
Applying Borell-TIS inequality (see, e.g., \cite{20lectures}) we have
as $u \to \IF$
\bqny{
\pk{\int\limits_{[0,\IF)\backslash [t_1-\ve,t_1+\ve]}
\ind(V_1(s)> u^{1-H})ds>T_u/u} &\le&
\pk{\exists t\in [0,\IF) \backslash[t_1-\ve,t_1+\ve]:
V_1(t)> u^{1-H}}
\\&\le&
e^{-\frac{(u^{1-H}-M)^2}{2m^2}},
}
where
$$M = \E{\sup\limits_{\exists t\in [0,\IF) \backslash[t_1-\ve,t_1+\ve]} V_1(t)}
< \IF, \quad
m^2 = \max\limits_{\exists t\in [0,\IF) \backslash[t_1-\ve,t_1+\ve]
 }\Var\{V_1(t)\}.$$
Since $\Var\{V_1(t)\}$ achieves its unique maxima at $t_1$ we obtain
by \eqref{ratio_uniform}
that
 $$e^{-\frac{(u^{1-H}-M)^2}{2m^2}}  = o(\pk{V_1(t_1)>u^{1-H}}),
\quad u \to \IF$$
and the claim follows from the asymptotics of $\psi_1(T_u,u)$
given in Proposition \ref{*prop_sojourn_app}.
\QED
\\

\textbf{Proof of \eqref{*infintcum}.}
To prove the claim it is enough to show that as $u \to \IF$ and then
$S \to \IF$
\bqny{
\pk{\int\limits_{[0,\IF)\backslash [ut_*-S,ut_*+S]}
\ind\big(B(t)-c_1t>q_1u,B(t)-c_2t>q_2u\big)dt>T } = o(\Cr(u)), \quad u \to \IF
.}
We have that the probability above does not exceed
\bqny{
\pk{\exists t \in [0,\IF)\backslash [ut_*-S,ut_*+S]:
B(t)-c_1t>q_1u,B(t)-c_2t>q_2u}.}
From the proof of Theorem 3.1 in \cite{Lanpeng2BM},
Case (3) and the final asymptotics of $\Cr(u)$ given in \eqref{*parisclaim2}
it follows that the expression above equals $o(\Cr(u))$, as $u \to \IF$ and then
$S \to \IF$.
\QED \\

\textbf{Proof of \eqref{*cumasym_elimin_u}.}
Define
$$G(u,x) = \pk{\int\limits_{-S}^{0}\ind\Big(Z_u(s)+\eta_*s>x\Big)ds+
\int\limits_0^S \ind\big(B_*(s)-c_1s>x\big)ds >T}.$$
First we show that
\bqn{\label{*case2AppBorellInt}
\int\limits_\R
G(u,x)e^{\frac{\eta x}{t_*}-\frac{x^2}{2ut_*}}dx
=
\int\limits_{-M}^M
G(u,x)e^{\frac{\eta x}{t_*}}dx+A_{M,u},
}
where $A_{M,u} \to 0$ as $u \to \IF$ and then $M \to \IF$.
We have
\bqny{ |A_{M,u}| &=&
|\int\limits_\R
G(u,x)e^{\frac{\eta x}{t_*}-\frac{x^2}{2ut_*}}dx
- \int\limits_{-M}^M G(u,x)e^{\frac{\eta x}{t_*}}dx|
\\&\le&
|\int\limits_{-M}^M G(u,x)
(e^{\frac{\eta x}{t_*}-\frac{x^2}{2ut_*}}-e^{\frac{\eta x}{t_*}})dx|
+ \int\limits_{|x|>M}G(u,x)e^{\frac{\eta x}{t_*}}dx
=: |I_1|+I_2.
}
Since the variance of $Z_u$ (see \eqref{*Z_u})
converges to those of BM we have by
Borell-TIS inequality
for $x>0$, large $u$ and some $C>0$
\bqn{\label{*G(u,x)_upper_bound}
G(u,x) &\le& \pk{\exists t \in [-S,0):(Z_u(t)+\eta_*t)>x}
+\pk{\exists t \in [0,S]:(B_*(t)-c_1t)>x}
\notag\\&\le&
\pk{\exists t \in [-S,0]:(Z_u(t)-\E{Z_u(t)})>x}
+\pk{\exists t \in [0,S]:B_*(t)>x}
\le
e^{-x^2/C}.
}
Let $u>M^4$. For $x\in [-M,M]$ it holds, that
$1-e^{-\frac{x^2}{2ut_*}} \le \frac{x^2}{2ut_*} \le \frac{1}{M}$ and
hence for $u>M^4$ by \eqref{*G(u,x)_upper_bound} we have as $M \to \IF$
\bqny{
|I_1| &\le&
\int\limits_{-M}^0e^{\frac{\eta x}{t_*}}(1-e^{-\frac{x^2}{2ut_*}})dx+
\int\limits_{0}^M
e^{-x^2/C+\frac{\eta x}{t_*}}(1-e^{-\frac{x^2}{2ut_*}})dx
\le \frac{1}{M}\Big(
\int\limits_{-\IF}^0 e^{\frac{\eta x}{t_*}}+
\int\limits_0^{\IF}e^{-x^2/C+\frac{\eta x}{t_*}}\Big)
\to 0
.}
For $I_2$ we have
\bqny{
I_2 \le
\int\limits_{-\IF}^{-M}e^{\frac{\eta x}{t_*}}dx+
\int\limits_M^\IF
e^{-x^2/C}
e^{\frac{\eta x}{t_*}}dx\to 0, \ \ \ M\to \IF,
}
hence
 \eqref{*case2AppBorellInt} holds.
Next we show that
\bqny{
G(u,x) & \to &
\pk{\int\limits_{-S}^{0}\ind\Big(B(s)+\eta_*s>x\Big)ds+
\int\limits_0^S \ind(B_*(s)-c_1s>x)ds >T}, \ \ \ u \to \IF \notag
}
that is equivalent with
\bqny{\lim\limits_{u \to \IF}
\pk{\int\limits_{-S}^{S}\ind\Big(X_u(s)>x\Big)ds>T}
=
\pk{\int\limits_{-S}^{S}\ind\Big(B(s)+k(s)>x\Big)ds>T}
, }
where  $k(s) = \ind(s<0)\eta_*s-\ind(s\ge 0)c_1s$ and
\bqny{
X_u(t) = (Z_u(t)+\eta_*t)\mathbb{I}(t<0)+(B_*(t)-c_1t)\mathbb{I}(t\ge 0).
}
We have for large $u$
\bqny{
\E{(X_u(t)-X_u(s))^2} =
\begin{cases}
|t-s|+|t-s|^2 & t,s\ge 0\\
-\frac{(s-t)^2}{ut_*}+|t-s|+\frac{x^2(t-s)^2}{u^2t_*^2}
-\frac{2x(t-s)^2\eta_*}{ut_*}+\eta_*^2(t-s)^2
 & t,s \le 0\\
|t-s|-\frac{s^2}{ut_*}+\frac{x^2s^2}{u^2t_*^2}
-\frac{2xs(\eta_*s+c_1t)}{ut_*}+(\eta_*s+c_1t)^2
 & s<0<t
\end{cases}
}
implying for all $u$ large enough, some
$C>0$ and $t,s \in [-S,S+T]$ that
\bqny{
\E{(X_u(t)-X_u(s))^2} \le C|t-s|.
}
Next, by Proposition 9.2.4 in \cite{20lectures}
the family $X_u(t), \ u>0,  \ t \in [-S,S+T]$ is tight in
$\mathcal{B}(C([-S,S+T]))$ (Borell $\sigma$-algebra in the space of
the continuous functions on $[-S,S+T]$ generated by the cylindric sets). \\
As follows from \eqref{*Z_u}, $Z_u(t)$ converges to $B(t)$ in
the sense of convergence finite-dimensional distributions as $u \to \IF$,
$t\in [-S,S+T]$.
Thus, by Theorems 4 and 5 in Chapter 5 in \cite{BylinskiiShiryaevBook}
the tightness and convergence of finite-dimensional distributions imply
weak convergence
$$X_u(t) \Rightarrow B(t)+k(t)=:W(t), \ \  t \in [-S,S+T].$$

By Theorem 11 (Skorohod), Chapter 5 in
\cite{BylinskiiShiryaevBook}
there exists a probability space $\Omega$,
where all random processes have the
same distributions, while weak convergence becomes convergence
almost sure. Thus, we assume that
$X_u(t)\to W(t)$ a.s. as $u \to \IF$
as elements of $C[-S,S]$ space with the uniform metric.
We prove that for all $x \in \R$
\bqn{\label{*cum_conv_as}
\pk{\limit{u}\int\limits_{-S}^S\ind(X_u(t)>x)dt =
\int\limits_{-S}^S\ind(W(t)>x)dt }=1
.}
Fix $x \in \R$.
We shall show that as $u \to \IF$ with probability 1
\bqn{\label{*leb_measure}
\mu_\Lambda\{t \in [-S,S]\!:\! X_u(\omega,t)\!>\!x\!>\!W(\omega,t)\}\!+\!
\mu_\Lambda\{t \in [-S,S]\!:\! W(\omega,t)\!>\!x\!>\!X_u(\omega,t)\} \to 0,
}
where $\mu_\Lambda$ is the Lebesgue measure.
Since for any fixed $\varepsilon>0$ for large $u$ and $t\in [-S,S]$
with probability one
$|W(t)-X_u(t)|<\ve$ we have that
\bqny{& \ &
\mu_\Lambda\{t \in [-S,S]: X_u(\omega,t)>x>W(\omega,t)\}
+
\mu_\Lambda\{t \in [-S,S]: W(\omega,t)>x>X_u(\omega,t)\}
\\&\le&
\mu_\Lambda\{t \in [-S,S]: W(\omega,t)\in [-\ve+x,\ve+x]\}.
}
Thus, \eqref{*leb_measure} holds if
\bqn{\label{*good_trajectory}
\pk{\lim\limits_{\varepsilon \to 0}
\mu_\Lambda\{t\in [-S,S]: W(t)\in [-\varepsilon+x,x+\varepsilon]\}=0}=1.
}
Consider the subset $\Omega_*\subset \Omega$ consisting of all $\omega_*$
such that
$$\lim\limits_{\ve \to 0}\mu_\Lambda\{t\in [-S,S]: W(\omega_*,t)
\in [-\varepsilon+x,x+\varepsilon]\}>0.$$
Then for each $\omega_*$ there exists the set
$\mathcal{A}(\omega_*) \subset [-S,S]$
such that $\mu_\Lambda\{\mathcal{A(\omega_*)}\}>0$
and for $t \in \mathcal{A(\omega_*)}$
it holds, that $W(\omega_*,t) = x$. Thus,
$$\pk{\Omega_*} = \pk{\mu_\Lambda\{t\in [-S,S]: W(t)=x\}>0},$$
the right side of the equation above
equals 0 by Lemma 4.2 in \cite{lemma_bms}.
Hence we conclude that \eqref{*good_trajectory} holds, consequently
\eqref{*leb_measure} and \eqref{*cum_conv_as} are true.
Since convergence almost sure implies convergence in distribution
we have by \eqref{*cum_conv_as} that
for any fixed $x\in \R$
\bqny{
\lim\limits_{u\to \IF}
\pk{\int\limits_{-S}^S\ind(X_u(t)>x)dt>T} =
 \pk{\int\limits_{-S}^S \ind(W(t)>x)dt>T}
.}
By the dominated convergence theorem we obtain
$$\int\limits_{-M}^M G(u,x)e^{\frac{\eta x}{t_*}}dx \to
\int\limits_{-M}^M
\pk{\int\limits_{-S}^{0}\ind\Big(B(s)+\eta_*s>x\Big)ds+
\int\limits_0^S \ind(B_*(s)-c_1s>x)ds >T}
e^{\frac{\eta x}{t_*}}dx, \ \ \ \ u \to \IF.$$
Thus, the claim follows from the line above and \eqref{*case2AppBorellInt}.
\QED
\newline

\textbf{Proof of \eqref{*cum_proof_main_int_H<1/2}.}
We have by the proof of Theorem 3.1 in \cite{Lanpeng2BM},
Case (3) and the final asymptotics of $\Cr(u)$ given in \eqref{*parisclaim2}
\bqny{
&\ &
\pk{\int\limits_{[0,\IF)\backslash [ut_*-u\delta_u,ut_*+u\delta_u]}
\ind(B_H(t)-c_1t>q_1u,B_H(t)-c_2t>q_2u)dt>T_u }
\\&\le&
\pk{\exists t \in [0,\IF)\backslash [ut_*-u\delta_u,ut_*+u\delta_u]:
B_H(t)-c_1t>q_1u,B_H(t)-c_2t>q_2u}
\\&=& o(\Cr(u)), \quad u \to \IF
}
and hence
\bqny{
\pk{\int\limits_{ [ut_*-u\delta_u,ut_*+u\delta_u]}
\ind(B_H(t)-c_1t>q_1u,B_H(t)-c_2t>q_2u)dt>T_u }
\sim \Cr(u), \quad u \to \IF.
}
The last probability above
is equivalent with $g_1(u)+g_2(u)$ as $u \to \IF$, this observation follows
from the application of the double-sum method, see the proofs
of Theorem 3.1, Case (3) $H<1/2$ in \cite{Lanpeng2BM} and Theorem 2.1
in \cite{SojournInfty} case i).
\QED
\\

\textbf{Proof of Proposition \ref{*prop_sojourn_app}.}
If $H=1/2$, then an equality takes place,
see \cite{SojournInfty}, Eq. $[5]$.
Assume from now on that $H\neq 1/2$.
First let \eqref{*assumption_T_u}
holds with $T>0$. We have for $c>0$ with
$\widetilde{M}(u) = u^{1-H}\frac{
c^H}{(1-H)^{1-H}H^H}$
(recall, $\ind_a(b) = \ind(b>a), \ a,b \inr$)
\bqny{
h_{T_u}(u) &:=&
\pk{\int\limits_ 0^\IF \ind(B_H(t)\!-\!ct>u)dt>T_u}
\\&=&
\pk{u (u^{\frac{1}{H}-2}\frac{c^2(1-H)^{2-\frac{1}{H}}
}{2^{\frac{1}{2H}}H^2}
)\int\limits_0^\IF \ind_{\widetilde{M}(u)}(
\frac{B_H(tu)\widetilde{M}(u)}{u(1+ct)} )dt>T\frac{c^2(1-H)^{2-\frac{1}{H}}
}{2^{\frac{1}{2H}}H^2}}
.}
Next we apply Theorem 3.1 in \cite{SojournInfty}
to calculate the asymptotics
of the last probability above as $u \to \IF$.
For the parameters in the notation therein we have
\bqny{& \ &
 \alpha_0=\alpha_\IF = H, \
\sigma(t) = t^H, \
\overleftarrow{\sigma}(t) = t^{\frac{1}{H}},
\ t^* = \frac{H}{c(1-H)},
\ A = \frac{c^H}{H^H(1-H)^{1-H}}, \
x =
T\frac{c^2(1-H)^{2-\frac{1}{H}}
}{2^{\frac{1}{2H}}H^2}
\\ &\ &
B = \frac{c^{2+H}(1-H)^{2+H}}{H^{H+1}}, \  M(u) = u^{1-H}\frac{
c^H}{(1-H)^{1-H}H^H}, \
v(u) = u^{\frac{1}{H}-2}\frac{c^2(1-H)^{2-\frac{1}{H}}
}{2^{\frac{1}{2H}}H^2}.
}
and hence we obtain
\bqn{\label{*h_T_u}h_{T_u}(u)
\sim K_H\mathcal{B}_{2H}(TD)(C_H u^{1-H})^{\frac{1}{H}-1}
\Psi(C_Hu^{1-H}), \ \ \ u \to \IF,
}
 where
$$C_H = \frac{c^H}{H^H(1-H)^{1-H}} \quad
\text{ and } \quad D = 2^{-\frac{1}{2H}}c^2H^{-2}
(1-H)^{2-1/H}.$$

Assume that \eqref{*assumption_T_u} holds with $T = 0$.
For $\ve>0$ for all large $u$ we have $
h_{\ve u^{1/H-2}}(u)\le h_{T_u}(u)\le h_0(u)$ and thus
\bqny{
K_H\mathcal{B}_{2H}(\ve D)(C_H u^{1-H})^{\frac{1}{H}-1}
\Psi(C_Hu^{1-H}) \le h_{T_u}(u) \le
K_H\mathcal{B}_{2H}(0)(C_H u^{1-H})^{\frac{1}{H}-1}
.}
Since $\mathcal{B}_{2H}(\cdot)$ is a continuous function (Lemma 4.1 in
\cite{SojournInfty}) letting $\ve \to 0$ we obtain
\eqref{*h_T_u} for any $T_u$ satisfying \eqref{*assumption_T_u}.
Replacing in \eqref{*h_T_u} $u$ and $c$ by $q_1u$ and $c_1$ we obtain
the claim.
\QED

\bibliography{sojourn_2dim_archiv}{}
\bibliographystyle{apalike}

\end{document}